\input amstex
\documentstyle{amsppt}
\topmatter
\magnification=\magstep1
\pagewidth{5.2 in}
\pageheight{6.7 in}
\abovedisplayskip=12pt \belowdisplayskip=12pt
\NoBlackBoxes
\title An invariant $p$-adic $q$-integral on $\Bbb Z_p$\endtitle

\author Taekyun Kim\endauthor
\affil{ {\it EECS, Kyungpook National University,
 Taegu, 702-701, S. Korea\\
 e-mail: tkim64$\@$hanmail.net}}
\endaffil

\abstract{ In this paper, we study $p$-adic q-integral on $\Bbb
Z_p$ and give the integral equations related to $p$-adic
$q$-integral. From these integral equations, we derive some
interesting Witt's formulae related to Bernoulli and Euler
numbers.
 } \endabstract
 \subjclass  11S80 \endsubjclass
\thanks  \endthanks
\leftheadtext{     } \rightheadtext{ }
\endtopmatter

\document

\head 1. Introduction \endhead Let $p$ be a fixed odd prime.
Throughout $\Bbb Z_p$, $\Bbb Q_p$ and $\Bbb C_p$ will respectively
denote the rings of $p$-adic integers, the fields of $p$-adic
numbers and the completion of  of the algebraic closure of $\Bbb
Q_p$. For $d$ a fixed positive integer with $(p,d)=1$, let
$$X=X_d=\varprojlim_N \Bbb Z/dp^N , \;\;X_1=\Bbb Z_p,$$
$$X^*=\bigcup\Sb 0<a<dp\\ (a,p)=1\endSb a+dp\Bbb Z_p,$$
$$a+dp^N\Bbb Z_p=\{x\in X\mid x\equiv a\pmod{dp^N}\},$$
where $a\in \Bbb Z$ lies in $0\leq a<dp^N ,$ (cf. [1], [2]).

The $p$-adic absolute value in $\Bbb C_p$ is normalized so that
$|p|_p=\frac1p .$ Let $q$ be variously considered as an
indeterminate a complex number $q \in \Bbb C$, or a $p$-adic
number $q\in\Bbb C_p .$  If $q\in\Bbb C_p,$ we assume $|q-1|_p <
1$ for $|x|_p \le 1.$ Throughout this paper, we use the following
notation :
$$[x]_q=[x:q]= \frac{1-q^x}{1-q}=1+q+q^2+\cdots+q^{x-1}.$$
We say that $f$ is uniformly differentiable function at a point
 $ a \in \Bbb Z_p$-- and denote this property by $ f \in UD(\Bbb
 Z_p )$-- if the difference quotients
 $$ F_f (x, y)=\frac{f(x)-f(y)}{x-y} ,$$ have a limit $l=
 f^{\prime}(a)$ as $( x, y) \rightarrow (a,a) ,$ cf. [1, 2].
For $f\in U D(\Bbb Z_p),$ let us start with the expression
$$\frac1{[p^N]_q} \sum_{0 \le j < p^N} q^i f(j) =\sum_{0\le j<p^N} f(j)
\mu_q(j+p^N\Bbb Z_p), \text{ cf. [1, 2], }$$ representing
$q$-analogue of Riemann sums for $f$.

The integral of $f$ on $\Bbb Z_p$ will be defined as limit
($n\rightarrow \infty$) of these sums, when it exists. An
invariant $p$-adic $q$-integral of a function $f \in {\text{UD}}
(\Bbb Z_p) $ on $\Bbb Z_p$ is defined by
$$I_{q}(f)=\int_{X}f(x)d\mu_q(x)=\int_{\Bbb Z_p} f(x) d\mu_q (x)= \lim_{N\rightarrow \infty} \frac1{[p^N]_q}
\sum_{0\le j<p^N} f(j)q^j, \text{ see [2] }.$$ Note that if $f_n
\rightarrow f$ in $UD( Z_p)$; then
$$\int_{\Bbb Z_p} f_n(x) d \mu_q(x) \rightarrow \int_{\Bbb Z_p} f(x) d
\mu_q(x).$$
 The purpose of this paper is to give the integral equations
 related to $I_{q}(f)$ and to investigate some properties for
 $I_{-q}(f)$. From these integral equations, we derive the
 interesting formulae related to Bernoulli and Euler numbers.

 \head \S 2. $p$-adic invariant integrals on $\Bbb Z_p$ \endhead

For $f\in UD(\Bbb Z_p),$ the $p$-adic $q$-integral was defined by
$$I_q(f)=\lim_{N\rightarrow \infty}\frac{1}{[p^N]_q}\sum_{0\leq
x<p^N}q^xf(x)=\lim_{N\rightarrow \infty}\sum_{0\leq
x<p^N}f(x)\mu_q(x+p^N\Bbb Z_p)=\int_{\Bbb Z_p}f(x)d\mu_q(x),\tag 1
$$ representing $p$-adic $q$-analogue of Riemann integral for $f$ (see
[2]). Let $f_1(x)$ be translation with $f_1(x)=f(x+1)$. Then we
have
$$I_q(f_1)=\frac{1}{q}I_q(f)+\left(\frac{q-1}{\log
q}f^{\prime}(0)+(q-1)f(0)\right).\tag2$$
From Eq.(1), we drive
$$ q^nI_q(f_n)=I_q(f)+\frac{q-1}{\log
q}\left(\sum_{i=0}^{n-1}f^{\prime}(i)q^i +\log q
\sum_{i=0}^{n-1}f(i)q^i\right), \text{ where $n\in\Bbb N $,
$f_n(x)=f(x+n)$. }$$
  Therefore we obtain the following theorem:
\proclaim{ Theorem 1} For $f\in UD(\Bbb Z_p),$ $n\in\Bbb N,$ let
$f_n(x)=f(x+n)$. Then we have
$$q^nI_{q}(f_n)-I_{q}(f)=\frac{q-1}{\log
q}\left(\sum_{i=0}^{n-1}f^{\prime}(i)q^i +\log q
\sum_{i=0}^{n-1}f(i)q^i\right) $$
\endproclaim
Let $f(x)=e^{tx}$. From Theorem 1, we can derive the following
$$I_q(e^{tx})=\int_{\Bbb Z_p}e^{tx}d\mu_q(x)=\left(\frac{\log q+t}{qe^t-1}\right)\frac{q-1}{\log
q}.\tag3$$ In [1,3], the q-Bernoulli numbers was defined by
$$\frac{\log q+t}{qe^t-1}=\sum_{n=0}^{\infty}B_{n, q}\frac{t^n}{n!}.
\tag4$$ Let $\chi$ be Dirichlet's character with conductor
$d\in\Bbb N$ and let $f(x)=e^{tx}\chi(x)$. By Theorem 1, we see
that
$$\int_{\Bbb
Z_p}\chi(x)e^{xt}d\mu_{q}(x)=\frac{\sum_{a=0}^{d-1}(\chi(a)e^{ta}q^a
t+\log q \chi(a)q^ae^{ta})}{q^de^{dt}-1}\frac{q-1}{\log q}.
\tag5$$ In [3], the generalized $q$-Bernoulli numbers attached to
$\chi$ as follows:
$$\frac{\sum_{a=0}^{d-1}(\chi(a)e^{ta}q^a
t+\chi(a)q^a\log
qe^{ta})}{q^de^{dt}-1}=\sum_{n=0}^{\infty}B_{n,q,\chi}\frac{t^n}{n!}.
\tag6$$ By (3), (4) and (5), we obtain the following:
 \proclaim{ Theorem 2}
For $n \in \Bbb N$, we have
$$\int_{\Bbb Z_p}\chi(x)x^nd\mu_q(x)=\frac{q-1}{\log
q}B_{n,q,\chi}, \text{ and } \int_{\Bbb
Z_p}x^nd\mu_q(x)=\frac{q-1}{\log q}B_{n,q}.$$ \endproclaim Let us
consider $p$-adic $q$-integral in the sense of fermionic  as
follows:
$$I_{-q}(f)=\lim_{n\rightarrow
\infty}\frac{1}{[p^N]_{-q}}\sum_{x=0}^{p^N-1}f(x)(-q)^x=\int_{\Bbb
Z_p}f(x)d\mu_{-q}(x). \tag7$$ Let $f_1(x)$ be translation with
$f_1(x)=f(x+1).$ From (7), we derive
$$qI_{-q}(f_1)=-I_{-q}(f)+[2]_qf(0). \tag8$$
By (7) and (8), we easily see that
$$q^{n}I_{-q}(f_n)=(-1)^nI_{-q}(f)+[2]_q\sum_{l=0}^{n-1}(-1)^{n-1-l}q^lf(l),
\text{ where $n\in\Bbb N$, $f_n(x)=f(x+n)$.}$$ Therefore we obtain
the following:
 \proclaim{ Theorem 3} For $ f\in UD(\Bbb Z_p)$, we
have
$$q^{n}I_{-q}(f_n)=(-1)^nI_{-q}(f)+[2]_q\sum_{l=0}^{n-1}(-1)^{n-1-l}q^lf(l),$$
where $n\in\Bbb N$, $f_n(x)=f(x+n)$.

In particular, if $n$ is odd positive integer, then we have
$$q^{n}I_{-q}(f_n)+I_{-q}(f)=[2]_q\sum_{l=0}^{n-1}(-1)^{l}q^lf(l).$$
\endproclaim
Let $\chi$ be Dirichlet's character with conductor $d(=odd)\in
\Bbb N$ and let $f(x)=e^{tx}\chi(x)$. By Theorem 3, we easily see
that
$$\int_{\Bbb
Z_p}e^{tx}\chi(x)d\mu_{-q}(x)=\frac{[2]_q\sum_{l=0}^{d-1}(-1)^lq^l\chi(l)e^{lt}}{q^de^{dt}+1}=\sum_{n=0}^{\infty}
H_{n,\chi}(-q^{-1})\frac{t^n}{n!},$$ where $H_{n,\chi}(q)$ are
called generalized Frobenius-Euler numbers attached to $\chi$,
cf.[4,5].

From (8), we can also derive
$$I_{-q}(e^{tx})=\int_{\Bbb
Z_p}e^{tx}d\mu_{-q}(x)=\frac{1+q}{qe^t+1}=\sum_{n=0}^{\infty}H_n(-q^{-1})\frac{t^n}{n!},$$
where $H_n(q)$ are called Frobenius-Euler numbers, cf.[5,6,7,8,9].

\enddocument